\input amstex
\mag=\magstep1
\documentstyle{amsppt}
\nologo
\NoBlackBoxes
\NoRunningHeads

\def\fp{\flushpar}
\define\bk#1{\text{\tenptmbit #1}}
\define\eightptbk#1{\text{\eightptmbit{#1}}}
\define\inbox#1{$\boxed{\text{#1}}$}
\define\underbarl#1{\lower 1.4pt \hbox{\underbar{\raise 1.4pt \hbox{#1}}}}
\define\llg{\langle\hskip -1.5pt\langle}
\define\rrg{\rangle\hskip -1.5pt\rangle}

\define\lr#1{^{\sssize\left(#1\right)}}
\define\br#1{^{\sssize\left[#1\right]}}
\define\do#1{^{\sssize\left<#1\right>}}
\define\ord{\text{ord}}
\define\twolower#1{\hbox{\raise -2pt \hbox{$#1$}}}

\define\hookdownarrow{\kern 0pt
 \hbox{\vbox{\offinterlineskip \kern 0pt
  \hbox{$\cap$}\kern 0pt
  \hbox{\hskip 3.5pt$\downarrow$}\kern 0pt}}
}  %

\font\sixptrm=cmr6

\font\twelveptbf=cmbx12 

\def\boxit#1{\vbox{\hrule\hbox{\vrule\kern3pt
    \vbox to 43pt{\hsize 182pt\kern3pt#1\eject\kern3pt\vfill}
    \kern1pt\vrule}\hrule}} %
\def\oVrule{\vrule width .5pt}
\def\oHrule{\hrule height .5pt}
\def\lbox#1#2#3{\kern 0pt %
\dimen1=#1pt \dimen2=#2pt%
\advance\dimen1 by -1.0pt%
\dimen3=\dimen1%
\advance\dimen3 by -16pt%
\advance\dimen2 by -1.0pt%
\dimen4=\dimen2%
\advance\dimen4 by -20pt%
\vbox to #1pt
{
\hsize #2pt
\oHrule
\hbox to #2pt
{
\vsize \dimen1
\oVrule
\vbox to \dimen1
{
\hsize \dimen2
\vskip 8pt
\hbox to \dimen2
{
\vsize \dimen3
\hskip 10pt
%
%
\vbox to \dimen3{\hsize \dimen4\fp #3\vfil}%
\hfil\hskip 10pt
}%
\vskip 8pt
}%
\oVrule%
}%
\oHrule%
}\kern 0pt%
}%
\def\cbox#1#2#3{\kern 0pt %
\dimen1=#1pt \dimen2=#2pt%
\advance\dimen1 by -1.0pt%
\dimen3=\dimen1%
\advance\dimen3 by -3pt%
\advance\dimen2 by -1.0pt%
\dimen4=\dimen2%
\advance\dimen4 by -3pt%
\vbox to #1pt%
{%
\hsize #2pt
\oHrule%
\hbox to #2pt%
{%
\vsize \dimen1%
\oVrule%
\vbox to \dimen1%
{%
\hsize \dimen2%
\vskip 1.5pt
\hbox to \dimen2%
{%
\vsize \dimen3%
\hskip 1.5pt\hfil
\vbox to \dimen3{\hsize \dimen4\vfil\hbox{#3}\vfil}%
\hfil\hskip 1.5pt
}%
\vskip 1.5pt
}%
\oVrule%
}%
\oHrule%
}\kern 0pt%
}%
\def\cboxit#1#2#3{$\hbox{\lower 2.5pt \hbox{\cbox{#1}{#2}{#3}}}$}
\document
\baselineskip=13.0pt 
{}
\rightline{\eightpoint 22nd December, 2003}
\vskip 10pt
\centerline{\twelveptbf Kummer's Original Type Congruence Relation}
\vskip 3pt
\centerline{\twelveptbf for the Universal Bernoulli Numbers}
\vskip 20pt
\centerline{Yoshihiro \^ONISHI}
\vskip 25pt
\fp
{\bf Preface.}\ \ \ 
\vskip 5pt
\fp
The aim of this paper is 
to give a congruence on universal Bernoulli numbers 
which congruence is the same type of Kummer's original paper \cite{{\bf K}}.  
The remarkable thing is the index of prime power of the modulus of 
the congruence is the half of the original one.  
We mention in this paper that this estimate is best possible.  
It is surprising fact for the author 
that the critical index is {\it not less than half} of the original.  

The motivation of this work is the investigation on 
generalized Bernoulli-Hurwitz numbers by the author himself in \cite{{\bf \^O}}.  
Kummer's original type congruences in \cite{{\bf \^O}} hold 
modulo the same power as the original one.  
When the author was working to get a proof of such 
Kummer type congruences in \cite{{\bf \^O}}, he knew several researches 
on universal Bernoulli numbers, 
especially Adelberg's remarkable papers \cite{{\bf A1}}  and  \cite{{\bf A2}}.  

There are three key lemmas for our proof of the main result,  
namely Lemma 3.2.1, Lemma 3.2.8, and Lemma 3.3.1.  
Lemma 3.2.8 is a very natural extension of already proved Lemma 3.2.1.  
Reglettablly, Lemma 3.2.8 is not yet proved.  
The Lemma 3.2.8 might be true by several reason.  
The author hope that 3.2.8 would be proved in the near future.  

Professor Adelberg gave several crucial comments on the 
first version of this paper.  
This version is much improved by his comments.  
The number theorists in Japan did not know  
the univesal Bernoulli numbers.  
The author expects that this research would be useful 
for the reserchers who are interested in 
Bernoulli numbers and Hurwitz numbers.  
\vskip 13pt

\newpage 

\fp
{\bf Contents}
\fp
\vskip 7pt
\fp
\inbox{1}   Combinatorics\fp\ \ \ 
 1.1. Properties of factorial\fp\ \ \ 
 1.2. Lagrange inversion formula
\vskip 3pt
\fp
\inbox{2}   Universal Bernoulli numbers and thier properties\fp\ \ \ 
 2.1. Definition of universal Bernoulli numbers\fp\ \ \ 
 2.2. Schur function type expression of universal Bernoulli numbers\fp\ \ \ 
 2.3. Clarke's theorem
\vskip 3pt
\fp
\inbox{3}  Kummer type congruence for universal Bernoulli numbers\fp\ \ \ 
 3.1. Main results\fp\ \ \ 
 3.2. Preparation for the proof of main theorem (1)\fp\ \ \ 
 3.3. Preparation for the proof of main theorem (2)\fp\ \ \ 
 3.4. Proof of the main theorem\fp\ \ \ 
 3.5. Kummer-Adelberg congruence
\vskip 3pt
\fp
References
\vskip 13pt


\fp
{\bf Convention}. 
\vskip 5pt
\fp
(1) For a rational number  $\alpha$, we denote by  $\lfloor \alpha \rfloor$  
the greatest integer not exceed  $\alpha$,  
and denote by  $\lceil \alpha \rceil$  
the least integer not less than  $\alpha$.  
\vskip 5pt
\fp
(2) We use the notation
  $$
  (n)_r=n(n-1)\cdots(n-r+1).  
  $$
\vskip 5pt
\fp
(3) The congruence relations on polynomials in several variables 
means the congruence on the coefficients of each the similar terms.  
\vskip 5pt
\fp
(4) Let  $R$  be a commutative ring and  $z$  be an indeterminate.  
We denote 
  $$
  R\llg z \rrg=\bigg\{\sum_{n=0}^{\infty}a_n\frac{z^n}{n!}\ \bigg|\ 
                    a_n\in R\ \ \text{for all  $n$}\bigg\}. 
  $$
For two power series 
  $$
  \varphi(z)=\sum_{n=0}^{\infty}a_n\frac{z^n}{n!}
  \ \ \ \ \text{and}\ \ \ \ 
  \psi(z)=\sum_{n=0}^{\infty}b_n\frac{z^n}{n!}
  $$
in  $R\llg z \rrg$  and for a subset  $S\subset R$, the congruence 
  $$
  \varphi(z)\equiv \psi(z) \mod{S}
  $$
means that  $a_n-b_n\in S$  for all  $n$.  
Especially, if  $p$  is a prime element in  $R$, 
and  $S=p^dR$  for a positive integer  $d$, 
then we write simply that  $\varphi(z)\equiv \psi(z)$ mod  $p^d$.  
\fp

\newpage 

\fp
\inbox{\bf 1}\ \ 
{\bf Combinatorics. } \ \ 
\vskip 5pt
\fp
{\bf 1.1. \  Properties of factorial. }\ \     
We recall the following two properties about factorial.  
Let  $n$  and  $k$  be non-negative integers, 
and  $p$  be a rational prime.  
If  $n=kp+a$  ($0\leqq a <p$),  then
  $$
  \ord_p(n!)=\ord_p((kp)!)=\ord_p(k!)+k. 
  \tag{1.1.1}
  $$
We denote by  $S_p(n)$  the sum of $p$-adic digits of  $n$.  
It is well known that 
  $$
  \ord_p(n!)=\frac{n-S_p(n)}{p-1}. 
  \tag{1.1.2}
  $$
\vskip 10pt

\fp
{\bf 1.2. \  Lagrange inversion formula. }\ \     
For a power series  $F(z)$, 
we denote by  $[z^n]F(z)$  the coefficient of the term  $z^n$.  
The following is called Lagrange inversion formula.  
This is a very strong tool.  
\vskip 5pt
\fp
\lbox{85}{360}
{{\bf Proposition 1.2.1.}\ \ 
Let  $\varphi(u)=u+\cdots$  is a power series with only positive terms in  $u$.  
The coefficient of degree  $1$  term is supposed to be  $1$.  
Let  $\psi(t)=\varphi^{-1}(t)$  be the inverse power series of  $\varphi(u)$,  
namely  $\varphi(\psi(t))=t$.  Then
  $$
  [u^n]\Big(\frac{u}{\varphi(u)}\Big)^{\ell}
  =\frac{\ell}{\ell-n}[t^n]\Big(\frac{\psi(t)}{t}\Big)^{\ell-n}. 
  $$
}
\vskip 5pt
\fp
About the proof of this, 
see the reference in  \cite{{\bf A1}}，p.123，Proposition 2.1.  
\vskip 13pt


\fp
\inbox{\bf 2}\ \ 
{\bf The universal Bernoulli numbers and their properties. }\ \ 
\vskip 5pt
\fp
{\bf 2.1. Definition of the universal Bernoulli numbers. }\ \     
Let  $c_1$, $c_2$, $\cdots$  are indeterminates.  
We cosider the power series
  $$
  u=u(t)=t+\sum_{n=1}^{\infty}c_{\sssize n}\frac{t^{n+1}}{n+1}
  \tag{2.1.1}
  $$
and its inverse power series
  $$
  t=t(u)=u-c_1\frac{u^2}{2!}+(3{c_1}^2-2c_2)\frac{u^3}{3!}+\cdots. 
  \tag{2.1.2}
  $$
We define the {\it universal Bernoulli numbers} (of degree $1$)   
$\hat{B}_n\in \bold Q[c_1,c_2,\cdots]$  by  
  $$
  \frac{u}{t(u)}
  =\sum_{n=0}^{\infty}\hat{B}_n\frac{u^n}{n!}
  \tag{2.1.2}
  $$

If we specialize  $c_n$  as  $c_n=(-1)^n$,  
we have  $\hat{B}_n=B_n$, the usiual Bernoulli numbers,  
because of  $\psi(t)=\log(1+t)$  and  $\varphi(u)=e^t-1$.  
\vskip 13pt


\fp
{\bf 2.2. \ Schur function type expression of universal Bernoulli numbers. }\ \   
\fp
We introduce several notations.  
For a finite sequence  $U=(U_1,U_2,\cdots)$  of non-negative integers, 
we define the {\it weight} of  $U$  to be  $w(U)=\sum_jjU_j$,  
and the {\it degree} of  $U$  to be  $d(U)=\sum_jU_j$.  
We can regard  $U$  to be a partition of  $w(U)$.  
For simplicity, we use notations
  $$
  \aligned
  U!&=U_1!U_2!\cdots\ \ \ \binom{d}{U}=\frac{d!}{U!},  \\
  \Lambda^{\sssize U}&=2^{\sssize U_1}3^{\sssize U_2}4^{\sssize U_3}\cdots,  \\
  c^{\sssize U}&={c_1}^{\sssize U_1}{c_2}^{\sssize U_2}{c_3}^{\sssize U_3}\cdots. 
  \endaligned  
  \tag{2.2.1}
  $$
Moreover we denote 
  $$
  \gamma_{\sssize U}
  =\Lambda^{\sssize U}U!.
  \tag{2.2.2}
  $$
For the function  $\psi(t)$  in 1.2.1, 
let  $h(t)=(\psi(t)/t)-1$.  Then
  $$
  (\psi(t)/t)^s
  =(1+h(t))^s
  =\sum_{d=0}^{\infty}\binom{s}{d}h^d(t), 
  \tag{2.2.3}
  $$
and 
  $$
  h^d(t)=\sum_{d(U)=d}\binom{d}{U}
  \frac{c^{\sssize U}}{\Lambda^{\sssize U}}\hskip 2pt t^{\sssize w(U)}.  
  \tag{2.2.4}
  $$
Therefore, by denoting 
  $$
  \tau_{\sssize U}=(-1)^{d(U)-1}\frac{(w(U)+d(U)-2)!}{\gamma_{\sssize U}}. 
  \tag{2.2.5}
  $$
and by using  Proposition 1.2.1 for  $\ell=1$, we have
\vskip 5pt
\fp
\lbox{62}{360}
{{\bf Proposition 2.2.6.}\ \ 
  $$
  \frac{\hat{B}_n}{n}
   =\hskip -2pt\sum_{w(U)=n}
   \tau_{\sssize U}c^{\sssize U}．
  $$
}
\vskip 5pt
\fp
According to Haigh's pointing out (\cite{{\bf C}}, p.594, $\ell$.$-$12), 
we call this expression {\it Schur function expression} of  $\hat{B}_n$.  
\vskip 20pt


\fp
{\bf 2.3. Clarke's theorem.}
\fp
For convenience of the reader, 
we mention here Clarke's theorem that is 
the universal version of von Staudt-Clausen theorem 
jointed with von Staudt second theorem.  

We denote  $a|_p=a/p^{\ord_pa}$  for given positive integer  $a$.  
\vskip 5pt
\lbox{223}{360}
{{\bf Proposition 2.3.1.}\ \ 
  $$
  \align
  \hat{B}_1&=\frac12c_1, \\
  \frac{\hat{B}_2}{2}&=-\frac14{c_1}^2+\frac13c_2,  \\
  \frac{\hat{B}_n}{n}
  &\equiv 
  \cases
  \dsize\sum_{\Sb n=a(p-1) \\ p:\text{\sixptrm prime}\endSb}\hskip -10pt 
  \dfrac{{a|_{\sssize p}}^{-1}\ \text{mod}\ p^{1+\ord_pa}}{p^{1+\ord_pa}}\ 
  {c_{p-1}}^a
      \hskip 30pt  (\text{if  $n\equiv 0$ mod $4$}) & \\
  \dfrac{{c_1}^{n-6}{c_3}^2}{2}-\dfrac{n{c_1}^n}{8}
  + \hskip -10pt
  \dsize\sum_{\Sb n=a(p-1) \\ p:\text{\sixptrm odd prime}\endSb}\hskip -10pt 
  \dfrac{{a|_{\sssize p}}^{-1}\ \text{mod}\ p^{1+\ord_pa}}{p^{1+\ord_pa}}\ 
  {c_{p-1}}^a & \\
      \hskip 123pt 
      (\text{if  $n\neq 2$  and  $n\equiv 2$  mod $4$}) & \\
  \dfrac{{c_1}^n+{c_1}^{n-3}c_3}{2} \hskip 51pt
      (\text{if  $n\neq 1$  and  $n\equiv 1$, $3$ mod  $4$}) & 
  \endcases  \\
  &\mod \bold Z[c_1, c_2, \cdots]. 
  \endalign
  $$
}
\fp
\vskip 5pt
\fp
The proof is given by analysing the Schur function expression 2.2.6
of  $\hat{B}_n$. 
See \cite{{\bf C}}, Theorem 5.  
\vskip 18pt


\fp
\inbox{\bf 3}\ \ 
{\bf Kummer type congruence relations for universal Bernoulli numbers. }\ \ 
\vskip 5pt
\fp
{\bf 3.1. Main result.}
\fp
We prove that Kummer's original type congruence 
for universal Bernoulli numbers holds 
modulo  $p^{\lfloor a/2 \rfloor}$  as follows.  
\vskip 5pt
\fp
\lbox{89}{360}
{{\bf Theorem 3.1.1.}\ \ 
Fix a prime number  $p$. 
Let  $a$  and  $n$  be positive integers such that 
$n>a$  and  $n\not\equiv 0$  mod  $(p-1)$.  
Then
  $$
  \sum_{r=0}^a\binom{a}{r}(-1)^r {c_{p-1}}^{a-r}
           \frac{\hat{B}_{n+r(p-1)}}{n+r(p-1)} 
  \equiv 0 \mod p^{\lfloor a/2 \rfloor}. 
  $$
}
\vskip 5pt
\fp
{\bf Remark 3.1.2.}\ \  
(1)  Suppose  $n>a=1$  and  $n\not\equiv 0$, $1$  mod  $(p-1)$.  
Then the above congruence holds modulo  $p$.  
This fact is proved by Adelberg (\cite{{\bf A1}}, Theorem 3.2).  
For the case of  $n\equiv 1$  mod  $(p-1)$,   
see the main theorem of \cite{{\bf A2}}.  \fp
(2)  Let  $p\geqq 7$  be a prime.  
Let  $U$  is the partition with  $U_1=p$,  $U_{2p-1}=(p-3)/2$, 
and  $U_j=0$  for the others.  
Then  $w(U)=p+(p-5)(2p-1)/2\equiv -1$  mod $(p-1)$.  
Then we can prove  
  $$
  \ord_p(\tau_{\sssize U})=(p-5)/2 \ \ 
  (=\lfloor (p-4)/2 \rfloor).
  $$  
Therefore we see that the above is a best possible estimate 
by the equation (3.4.1) below  for  $a=p-4$  and  $n=w(U)$.  
\fp
(3)  In the example in (2) above, if we set  $p=5$  
then we see  $\ord_5(\tau_{\sssize U})=0$.  
Here  $n=w(U)=5\equiv 1$  mod  $(5-1)$.  
So this is a case removed from Theorem 3.2 in \cite{{\bf A1}}.  
Keeping this example in mind, 
we can slightly improve Lemma 3.3.1 below. 
Then we can prove the case of  $a=1$  in 3.1.1.  
\vskip 6pt

We show in 3.5 that 3.1.1 implies directly Adelberg's congruence 
(\cite{{\bf A2}}, part (i) of Theorem).  
\vskip 3pt
\fp
\lbox{73}{360}
{{\bf Corollary 3.1.3.}\ (Adelberg's congruence relation)\   
Let  $n$  and  $a$  be positive integers.  
If  $n\not\equiv 0$, $1$  mod  $p-1$  and  $n>a$, then
  $$
  {c_{p-1}}^{p^{a-1}}\hskip -3pt\cdot
  \frac{\hat{B}_{n\vphantom{n+p^a(p-1)}}}{n}
  \equiv  
  \frac{\hat{B}_{n+p^{a-1}(p-1)}}{n+p^{a-1}(p-1)}
  \mod p^a. 
  $$
}
\vskip 19pt


\fp
{\bf 3.2. \ Preparation for the proof of main theorem (1). }\ \ \ 
\fp
From now on, we denote by  $\overline{k}$ 
the least non-negative integer that is congruent to  $k$  modulo $p$.   
In this section we prove the following estimate.  
\fp
\vskip 5pt
\fp
\lbox{112}{360}
{{\bf Lemma 3.2.1.}\ \ 
We fix an odd prime  $p$.  
Let  $a$, $q$, and  $n$  be non-negative integers.  Then 
  $$
  \sum_{r=0}^a\frac{\left((r+q)p+n\right)!}{(r+q)!p^{r+q}}
  \binom{a}{r}
  \equiv 0 \mod p^M, 
  $$
where
  $$
  M=\cases
    \ord_p(n!) \ \ & (\text{if  $n\geqq ap$}), \\
           a-\lfloor n/p \rfloor + \ord_p(n!)
          -  \lfloor(a-\lfloor n/p \rfloor)/p\rfloor 
                          &  (\text{if  $n<ap$}). 
    \endcases
  $$
}
\vskip 5pt
\fp
{\it Proof}.  We give a proof by using a generating function
\footnote{The author does not know any proofs of 3.2.1 
in combinatrics. }.  
We consider the function 
  $$
  F(v)=\exp(v^p/p) 
      =\sum_{r=0}^{\infty}\frac{(v^p/p)^r}{r!} 
      =\sum_{r=0}^{\infty}\frac{(rp)!}{r!p^r}\frac{v^{rp}}{(rp)!}. 
  \tag{3.2.2}
  $$
Obviously this series belongs to  $\bold Z\llg v \rrg$.  
We investigate 
  $$
  \Big(\Big(\frac{d}{dv}\Big)^p+1\Big)^a\big(v^nF(v)\big).  
  \tag{3.2.3}
  $$
This is a polynomial of  $v$  times  $F(v)$.  
The coefficients of the polynimial are as follows.  
The highest term is  $v^{n+ap(p-1)}$  with the coefficient  $1$, 
and
\fp
(a) \underbarl{if $n\geqq ap$}, then the lowest term is  $v^{n-ap}$,  and 
the whole polynomial is a polynomial of  $v^p$  times  $v^{n-ap}$; 
\fp
(b) \underbarl{if $n<ap$}, then the lowest term is  $v^{\overline{n}}$,  
the whole polynomial is a polynomial of  $v^p$  times  $v^{\overline{n}}$.  
\fp
We devide these terms into several groups as follows. 
The highewt term  $v^{n+ap(p-1)}$  is itself consists one of the groups. 
The following higher  $p$  terms, including the terms with zero coefficient, 
consist the next one of the groups.  
We continue similar grouping with  $p$  terms each.  
Although the number of the finally remaining terms is possibly 
less than  $p$, we regard the terms to be a group.  
If  $n\geqq ap$, then we devide 
  $$
  \frac{(n+ap(p-1))-(n-ap)}{p}+1=ap+1
  $$  
terms in to the groups.  Hence we get  $a+1$  groups.  
If  $n>ap$, then we devide 
  $$
  \align
  \{(n+ap(p-1))-\overline{n}\}/p
 =&\{(n+ap(p-1))-(n-\lfloor n/p \rfloor p)\}/p \\
 =&a(p-1)+\lfloor n/p \rfloor p
  \endalign
  $$
terms into the groups. Hence we get 
  $$
  \lceil (a(p-1)+\lfloor n/p \rfloor p)/p \rceil+1
 =a-\lfloor (a-\lfloor n/p \rfloor)/p \rfloor +1
  $$
groups. 
We denote by  $w_0(v)$, $w_1(v)$, $\cdots$  
the sums of the terms in the each group, 
according to the order from lower to higher. 
Then we have
  $$
  \Big(\Big(\frac{d}{dv}\Big)^p+1\Big)^a\big(v^nF(v)\big)
  =\cases
  \hskip 22pt
  \dsize\sum_{j=0}^{a}w_j(v) &  (\text{if  $n\geqq ap$}), \\
  \dsize\sum_{j=0}^{a-\lfloor(a-\lfloor \frac{n}{p} \rfloor)/p \rfloor}
     \hskip -20pt w_j(v) 
                         &  (\text{if  $n< ap$}). 
  \endcases
  \tag{3.2.4}
  $$
The most important things are as follows: namely
\fp
(a) \underbarl{if $n\geqq ap$}, then
  $$
  w_j(v)\hskip -1pt
  =\hskip -1pt\cases
  p^{\ord_p((n)_{\sssize (a-j)p})+a-j}v^{n-(a-j)p}\dsize\sum_{i=0}^{p-1}F_{ji}\hskip 1pt v^{ip}
       &  (\text{if  $0\leqq j<a$}) \\
  v^{n+ap(p-1)}
       &  (\text{if  $j=a$})
  \endcases
  \tag{3.2.5}
  $$
\fp
(b) \underbarl{if  $n< ap$}, then
  $$
  w_j(v)\hskip -1pt
  =\hskip -1pt\cases
  \hskip -2pt
  p^{\ord_p(n!)-\lfloor \frac{n}{p} \rfloor}
  p^{a-\lfloor(a-\lfloor \frac{n}{p} \rfloor)/p \rfloor)}
  v^{\overline{n}}\dsize\sum_{i=0}^{\overline{p-a}} F_{0i}\hskip 1pt v^{ip} 
           \hskip 10pt  (\text{if  $j=0$})& \\ 
  \hskip -2pt
  p^{\ord_p(n!)-\lfloor \frac{n}{p} \rfloor-j}
  p^{a-\lfloor(a-\lfloor \frac{n}{p} \rfloor)/p \rfloor)-j}
  v^{\overline{n}+jp^2+\overline{(p-a)}p}
                  \dsize\sum_{i=0}^{p-1} F_{ji}\hskip 1pt v^{ip} & \\
  \hskip -2pt
           \hskip 110pt  (\text{if  $1\leqq j < 
            \ord_p(n!)-\lfloor \frac{n}{p} \rfloor$}) & \\ 
  p^{a-\lfloor(a-\lfloor \frac{n}{p} \rfloor)/p \rfloor)-j}
  v^{\overline{n}+jp^2+\overline{(p-a)}p}
                  \dsize\sum_{i=0}^{p-1} F_{ji}\hskip 1pt v^{ip} & \\
  \hskip -2pt
           \hskip 37pt  (\text{if  $\ord_p(n!)-\lfloor \frac{n}{p} \rfloor
            \leqq j < a-\lfloor(a-\lfloor \frac{n}{p} \rfloor)/p\rfloor$}) 
                           &\\ 
  \hskip -2pt
  v^{n+ap(p-1)}
        \hskip 70pt
         (\text{if  $j=a-\lfloor(a-\lfloor \frac{n}{p} \rfloor)/p\rfloor$}). & 
  \endcases
  \tag{3.2.6}
  $$
While these facts can be proved by induction on  $a$, 
we minimally explain on the process.  
In the first two cases in  (3.2.6),  
the $p$-power of the head is appeared when we 
operate  $(\frac{d}{dv})^p$  to  $v^n$.  
Besides each oparation of  $(\frac{d}{dv})^p$  gives multiplication of one  $p$, 
this factor is spoiled by the contribution by the similar terms 
obtained by the oparation by  $(\frac{d}{dv})^p$  to  $F(v)$.  
Now we assume that  $\ell$  is a positive integer 
such that there exists a positive integer  $k$  with 
satisfying  $p^2|k$  and  $k\leqq \ell \leqq k+(p-1)$.   
Let us consider the following situation;
namely, after operating  $(\frac{d}{dv})^p$  several times,  
we are going to operating  $(\frac{d}{dv})^p$  to  $v^{\ell}$.  
If this operation is done the coefficients is multiplied 
by at least  $p^2$.  
In this situation, the similar terms come from 
the group of the next higher level.  
So the order of  $p$  is incresed by at least one.  
We just finish to explaine the first $p$-factors in the first 
two cases in  (3.2.6).  

The other  $p$-factors come by the following reason.  
When  $(\frac{d}{dv})^p+1$  operates to  $F(v)$, we have
  $$
  \aligned
  \Big(\Big(\frac{d}{dv}\Big)^p+1\Big)F(v)
  =&\Big(\frac{d}{dv}\Big)^{p-1}v^{p-1}F(v)+F(v) \\
  =&\big((p-1)!F(v)+\cdots\big)+F(v).
  \endaligned
  \tag{3.2.7}
  $$
Since  $(p-1)!+1\equiv 0$  mod  $p$, 
all the coefficients of the terms in (3.2.7) are divisible by  $p$.  
Only this mechanism gives rise to the other  $p$-factors.  

If we regard the right hand side of (3.2.4) to be 
a linear combination of terms  $\{v^m/m!\}$,   
the term whose coefficient has the least  $p$-factor is 
just the first terms in (3.2.5) and in (3.2.6).  
In other words, if we regard the right hand side to be 
a linear combination of  $\{v^jF(v)\}$,  
the term whose coefficient has the least  $p$-factor is 
the term  $v^jF(v)$  with the least $j$.  
Because the coefficient of  $v^{qp+n}/(qp+r)!$  in 
  $$
  \bigg(\Big(\frac{d}{dv}\Big)^p+1\bigg)^av^nF(v)
  =\sum_{r=0}^a\binom{a}{r}\Big(\frac{d}{dv}\Big)^{pr}
  \Big(\sum_{j=0}^{\infty}\frac{(jp+n)!}{j!p^r}
  \frac{v^{jp+n}}{(jp+n)!}\Big)
  $$
is just the left hand side of our claim,  
the proof has completed.  
\qed

\vskip 10pt
We need a variant with replacing  $q$  in 3.2.1  by negative  $r_0$  as follows. 
\vskip 3pt
\fp
\lbox{131}{360}
{{\bf Lemma 3.2.8.}\ \ (This is a Conjecture at present.)  
Let  $p$  be an odd prime, 
and  $a$  be a positive integer.  
Suppose  $r_0$  is an integer with  $0< r_0 \leqq a$.  
Let  $n\geqq r_0p$  be an integer.  
Then 
  $$
  \sum_{r=r_0}^a\frac{\left((r-r_0)p+n\right)!}{(r-r_0)!p^{r-r_0}}
  \binom{a}{r}
  \equiv 0 \mod p^M, 
  $$
where 
  $$
  M=\cases
    \ord_p(n!) \ \ &  (\text{if  $n\geqq ap$}), \\
           a-\lfloor n/p \rfloor + \ord_p(n!)
          -  \lfloor(a-\lfloor n/p \rfloor)/p\rfloor 
                          &  (\text{if $n<ap$}).  
    \endcases
  $$
}
\vskip 5pt
\fp
Althogh this lemma is not yet proved, many numerical examples 
suggest this would be true and 
it seems natural if we replase the factorials by the function  $\varGamma$  
with comparing 3.2.1.  
So we can strongly expect the truth of this lemma.  
\vskip 13pt


\fp
{\bf 3.3. Preparation for the proof of main theorem (2). }\ \ \ 
\fp
We show the following Lemma in this subsection.  
\vskip 3pt
\fp
\lbox{71}{360}
{{\bf Lemma 3.3.1.}\ \   
Let  $p$  be an odd prime, and  $U$  be a partition with  $U_{p-1}=0$.  
Assume  $d(U)>0$.  
Then for  $\tau_{\sssize U}$  defined in  (2.2.5) 
we have
  $$
  \ord_p (\tau_{\sssize U})
  \geqq \left\lfloor \frac{w(U)+d(U)-2}{2p} \right\rfloor. 
  $$
}
\vskip 5pt
\fp
{\bf Remark 3.3.2.} \ \ 
As is mentioned in 3.1.2(2), if  $U$  is the partition 
with  $U_1=p$, $U_{2p-1}=(p-5)/2$, and  $U_j=0$  for the others, 
then we have  $w(U)=p+\frac{(2p-1)(p-5)}{2}$,  $d(U)=p+\frac{p-5}{2}$, 
and  $\ord_p(\tau_{\sssize U})=(p-5)/2$.  
On the other hand 
  $$
  \bigg\lfloor \frac{w(U)+d(U)-2}{2p} \bigg\rfloor
  =\bigg\lfloor \frac{p^2-3p-2}{2p}\bigg\rfloor
  =\bigg\lfloor \frac{p-5}{2}+\frac{p-1}{p} \bigg\rfloor
  =\frac{p-5}{2}
  $$
Hence the above estimate is also best possible.  
\vskip 5pt
\fp
{\it Proof of {\rm 3.3.1}}. \ 
For simplicity we write  $w(U)=n$  and  $d(U)=d$.  
Since  $d(U)>0$, we have  $n+d-2>0$.  
Firstly we suppose  $U_{2p-1}\neq 0$.  
Then we have 
  $$
  \allowdisplaybreaks
  \align
 & \ord_p (\tau_{\sssize U})  
  =\ord_p ( (n+d-2)! )-\ord_p ( \gamma_{\sssize U} ) \\
 &=\ord_p \bigg( \Big(-2+\sum_{j\neq p-1}(j+1)U_j \Big)!\bigg)
    - \hskip -7pt
      \sum_{(\epsilon, k)\neq (1,1)}
      \hskip -10pt kU_{\epsilon p^k-1}
    - \sum_{j\neq p-1}\ord_p (U_j!) \\
 & \hskip 90pt 
   (\text{where  $\epsilon$  runs through the positive integers corprime to  $p$. }) \\
 & \geqq   
   \ord_p 
     \bigg(\Big(
      - 2
      + \hskip -7pt
        \sum_{j\neq p-1,~2p-1}\hskip -7pt jU_j 
      + 2p U_{2p-1}
     \Big)!\bigg)
      - \hskip -7pt
        \sum_{(\epsilon, k)\neq (1,1)}
        \hskip -7pt kU_{\epsilon p^k-1} 
      - \ord_p(U_{2p-1}!) \\
 &= \ord_p 
     \bigg(\Big(
      - 2
      + \sum_{p\not\hskip 2pt | j+1} jU_j 
      + \hskip -15pt
        \sum_{(\epsilon, k)\neq (1,1),~(2,1)}
        \hskip -15pt (\epsilon p^k-1)U_{\epsilon p^k-1}
      + 2p U_{2p-1}
     \Big)!\bigg) \\
 &  \hskip 50pt
   -\hskip -8pt
    \sum_{(\epsilon, k)\neq (1,1)}
    \hskip -8pt kU_{\epsilon p^k-1}
   -\ord_p(U_{2p-1}!) \\
 &\geqq  
    \sum_{\nu=1}^{\infty}
     \bigg\lfloor
      \frac{1}{p^{\nu}}
      \bigg(
      - 2 
      + \sum_{p\not\hskip 2pt | j+1} jU_j
      + \hskip -7pt
        \sum_{(\epsilon, k)\neq (1,1),~(2,1)} 
        \hskip -7pt (\epsilon p^k-1)U_{\epsilon p^k-1}
      + 2pU_{2p-1} 
      \bigg)
     \bigg\rfloor \\
 & \hskip 20pt
   -\hskip -7pt\sum_{(\epsilon, k)\neq (1,1)}
    \hskip -7pt kU_{\epsilon p^k-1} 
   -\ord_p(U_{2p-1}!) 
    \hskip 20pt (\because 
    \ord_p(N!)
    =\sum_{\nu=1}^{\infty}
     \big\lfloor \tfrac{N}{p^{\nu}} \big\rfloor）\\
 &=~
     \bigg\lfloor
      \frac{1}{p}
      \bigg(
      - 2 
      + \sum_{p\not\hskip 2pt | j+1} jU_j
      + \hskip -7pt
        \sum_{(\epsilon, k)\neq (1,1),~(2,1)} 
        \hskip -7pt (\epsilon p^k-1)U_{\epsilon p^k-1}
      + 2pU_{2p-1} 
      \bigg)
     \bigg\rfloor \\
 & \hskip 25pt
    +\sum_{\nu=2}^{\infty}
     \bigg\lfloor
      \frac{1}{p^{\nu}}
      \bigg(
      - 2 
      + \sum_{p\not\hskip 2pt | j+1} jU_j
      + \hskip -7pt
        \sum_{(\epsilon, k)\neq (1,1),~(2,1)} 
        \hskip -7pt (\epsilon p^k-1)U_{\epsilon p^k-1}
      + 2pU_{2p-1} 
      \bigg)
     \bigg\rfloor \\
 & \hskip 60pt 
    -\hskip -7pt\sum_{(\epsilon, k)\neq (1,1)}
     \hskip -7pt kU_{\epsilon p^k-1} 
    -\ord_p(U_{2p-1}!) ．
  \endalign
  $$
  $$
  \allowdisplaybreaks
  \align
 &\geqq 
     \bigg\lfloor
      \frac{1}{p}
      \bigg(
      - 2 
      + \sum_{p\not\hskip 2pt | j+1} jU_j
      + \hskip -7pt
        \sum_{(\epsilon, k)\neq (1,1),~(2,1)} 
        \hskip -7pt (\epsilon p^k-1)U_{\epsilon p^k-1}
      + 2pU_{2p-1} 
      \bigg)
     \bigg\rfloor \\
 &\hskip 25pt
    +\sum_{\nu=2}^{\infty}
     \bigg\lfloor
      \frac{- 2 + 2pU_{2p-1}}{p^{\nu}}
     \bigg\rfloor
    -\hskip -7pt\sum_{(\epsilon, k)\neq (1,1)}
     \hskip -7pt kU_{\epsilon p^k-1} 
    -\ord_p(U_{2p-1}!) \\
 &= \bigg\lfloor
      \frac{1}{p}
      \bigg(
      - 2 
      + \sum_{p\not\hskip 2pt | j+1} jU_j
      + \hskip -7pt
        \sum_{(\epsilon, k)\neq (1,1),~(2,1)} 
        \hskip -7pt (\epsilon p^k-1)U_{\epsilon p^k-1}
      + 2pU_{2p-1} 
      \bigg)
     \bigg\rfloor \\
   &\hskip 25pt
    -\hskip -7pt\sum_{(\epsilon, k)\neq (1,1)} 
     \hskip -7pt kU_{\epsilon p^k-1}
    +\sum_{\nu=2}^{\infty}
     \bigg\lfloor
      \frac{- 2 + 2pU_{2p-1}}{p^{\nu}}  
     \bigg\rfloor 
    -\ord_p(U_{2p-1}!) \\
 &= \bigg\lfloor
      \frac{1}{p}
      \bigg(
      - 2 
      + \sum_{p\not\hskip 2pt | j+1} jU_j
      + \hskip -7pt
        \sum_{(\epsilon, k)\neq (1,1),~(2,1)} 
        \hskip -7pt (\epsilon p^k-kp-1)U_{\epsilon p^k-1}
      + 2pU_{2p-1} 
      \bigg)
     \bigg\rfloor  \\
 &  \hskip 25pt
    -U_{2p-1}
    +\sum_{\nu=2}^{\infty}
     \bigg\lfloor
      \frac{- 2 + 2pU_{2p-1}}{p^{\nu}}  
     \bigg\rfloor 
    -\ord_p(U_{2p-1}!) \\
 &= \bigg\lfloor
      \frac{1}{p}
      \bigg(
      - 2 
      + \sum_{p\not\hskip 2pt | j+1} jU_j
      + \hskip -7pt
        \sum_{(\epsilon, k)\neq (1,1),~(2,1)} 
        \hskip -7pt (\epsilon p^k-kp-1)U_{\epsilon p^k-1}
      + 2pU_{2p-1} 
      \bigg)
     \bigg\rfloor \\ 
 &  \hskip 25pt
    - U_{2p-1}
    + \ord_p((-2+2pU_{2p-1})!)
    -\bigg\lfloor
      \frac{- 2 + 2pU_{2p-1}}{p}  
     \bigg\rfloor 
    - \ord_p(U_{2p-1}!) 
  \endalign
  $$
Now we can replace  $-\frac{2}{p}$  in the front by  $-\frac{1}{p}$.  
The reason is as follows:  
if the sum of the other terms in the parentheses just 
in  $\lfloor \ \ \rfloor$  is divisible by  $p$, 
then, after operating  $\lfloor \ \ \rfloor$,   
both of  $-\frac{2}{p}$  and  $-\frac{1}{p}$  give  $-1$;   
on the other hand, 
if such the sum of the terms is not divisible by  $p$,  
the remainder obtained by dividing it by  $p$  is 
at least  $\frac{1}{p}$.  
About the term   $jU_j$  in the second sum, we see  $j\geqq (j+1)/2$,  
and about the terms in the third sum, we see 
$\epsilon p^k-kp-1 > \epsilon p^k/2$  as far as  
$(\epsilon, k)\neq (1,1)$, $(2,1)$.  Therefore
  $$
  \allowdisplaybreaks
  \align
 &\geqq \bigg\lfloor
      \frac{1}{2p}
      \bigg(\hskip -3pt
      - 2 
      + \hskip -3pt\sum_{p\not\hskip 2pt | j+1} (j+1)U_j
      + \hskip -7pt
        \sum_{(\epsilon, k)\neq (1,1),~(2,1)} 
        \hskip -14pt \epsilon p^kU_{\epsilon p^k-1}
      + 2pU_{2p-1} 
      \bigg) 
      + U_{2p-1}
     \bigg\rfloor 
    -U_{2p-1} \\
 &  \hskip 25pt
    + \ord_p((-2+2pU_{2p-1})!)
    -\bigg\lfloor
      \frac{- 2 + 2pU_{2p-1}}{p}  
     \bigg\rfloor 
    - \ord_p(U_{2p-1}!) \\
 &= \bigg\lfloor
      \frac{-2+w(U)+d(U)}{2p}
     \bigg\rfloor \\
 &  \hskip 6pt
    +\ord_p(2U_{2p-1})!)+2U_{2p-1}-\ord_p(2pU_{2p-1})
    -\bigg\lfloor
      \frac{- 2 + 2pU_{2p-1}}{p}  
     \bigg\rfloor 
    - \ord_p(U_{2p-1}!) \\
&= \bigg\lfloor
      \frac{-2+w(U)+d(U)}{2p}
     \bigg\rfloor \\
 &  \hskip 6pt
    +\ord_p(2U_{2p-1})!)+2U_{2p-1}-\ord_p(2U_{2p-1})-1
    -(-1+2U_{2p-1}) 
    - \ord_p(U_{2p-1}!) \\
 &\geqq 
    \bigg\lfloor
      \frac{-2+w(U)+d(U)}{2p}
    \bigg\rfloor
   +\ord_p\bigg(\frac{(2U_{2p-1}-1)!}{U_{2p-1}!}\bigg)．
  \endalign
  $$
Hence we proved our statement if  $U_{2p-1}\neq 0$.   
If  $U_{2p-1}=0$,  we can get the same estimate by 
the similar argument by substituting  $U_{2p-1}=0$  
at the very beginning.  
\qed
\vskip 13pt


\fp
{\bf 3.4. Proof of the main theorem.}
\fp
We start now to prove the main result 3.1.1.  
By 2.2.6 (or \cite{{\bf C}}, p.594, Proposition 4), 
and substituting Schur function type expression of  $\hat{B}_n$, 
we have 
  $$
  \sum_{r=0}^a\binom{a}{r}(-1)^r {c_{p-1}}^{\hskip -8pt a-r}
           \frac{\hat{B}_{n+r(p-1)}}{n+r(p-1)} 
 =\sum_{r=0}^a\binom{a}{r}(-1)^r {c_{p-1}}^{\hskip -8pt a-r}
  \hskip -15pt\sum_{w(U)=n+r(p-1)}\hskip -10pt \tau_{\sssize U}c^U,    
  \tag{3.4.1}
  $$
where 
  $$
  \tau_{\sssize U}=(-1)^{d(U)-1}\frac{(w(U)+d(U)-2)!}{\gamma_{\sssize U}}. 
  $$
By enclosing as many as possible but less than or equal to  $r$  
factors  $c_{p-1}$  from  $c^U$, we have 
  $$
  \aligned
  =\sum_{r=0}^a\binom{a}{r}&(-1)^r{c_{p-1}}^{a-r} \\
  &\Bigg\{\hskip -2pt
    \sum_{w(U)=n}
         \hskip -10pt\tau_{\sssize U[r]}c^U{c_{p-1}}^r
   +\sum_{r_0=1}^r\hskip -15pt
         \hskip -10pt
           \sum_{\hskip 6pt\Sb\hskip 20pt w(U)=n+r_0(p-1) \\ U_{p-1}=0\endSb}
              \hskip -20pt
                 \tau_{\sssize U[r-r_0]}c^U{c_{p-1}}^{r-r_0}
  \Bigg\}, 
  \endaligned
  \tag{3.4.2}
  $$
where  $U[r]$  denotes the partition getting from  $U$  by adding  $r$  
to the  $(p-1)$-st entry.  
After exchanging the sum on  $r$  and the sum on  $U$, 
by writing down the terms  $\tau_{\sssize U[r]}$  and  $\tau_{\sssize U[r-r_0]}$, 
we see that 
  $$
  \aligned
 =&\sum_{w(U)=n}\hskip -5pt
         \frac{c^U{c_{p-1}}^{a}}{\gamma_{\sssize U|_{p-1}}} 
       \sum_{r=0}^a
        \binom{a}{r}(-1)^r (-1)^{d(U[r])+r-1} 
         \frac{\{w(U[r])+d(U[r])-2\}!}{p^{r+U_{p-1}}(r+U_{p-1})!} \\
   &\hskip 10pt
     + \ \ \sum_{r_0=1}^a\hskip -22pt
        \sum_{\Sb \hskip 24pt w(U)=n+r_0(p-1) \\ U_{p-1}=0 \endSb}
    \hskip -22pt
         \frac{c^U{c_{p-1}}^{a-r_0}}{\gamma_{\sssize U}} 
     \hskip 5pt\Bigg\{
         \sum_{r=r_0}^a\binom{a}{r}(-1)^r (-1)^{d(U[r-r_0])+r-1} \\
  &\hskip 130pt
        \cdot
        \frac{\{w(U[r-r_0])+d(U[r-r_0])-2\}!}{p^{r-r_0}(r-r_0)!} \Bigg\}, 
  \endaligned
  \tag{3.4.3}
  $$
where the symbol  $U|_{p-1}$  is  $U$  with out  $(p-1)$-entry;  
so that  $\gamma_{\sssize U|_{p-1}}$  written in the first sum 
means  $\gamma_{\sssize U}$  with neglected  the factors 
$p^{U_{p-1}}U_{p-1}!$  coming from its  $(p-1)$-st entry.  
We remark here that  
  $$
  \aligned
  \gamma_{\sssize U[r]|_{p-1}}
 &=\gamma_{\sssize U|_{p-1}},  \\
  \gamma_{\sssize U[r]|_{p-1}}p^{r+U_{p-1}}(r+U_{p-1})!
 &=\gamma_{\sssize U[r]}. 
  \endaligned
  \tag{3.4.4}
  $$
Note that 
  $$
  \gamma_{\sssize U}
  =\left(2^{U_1}\cdots(p-1)^{U_{p-2}}(p+1)^{U_p}\cdots\right)\cdot
               \left(U_1!\cdots U_{p-2}!U_p!\cdots\right)
  \tag{3.4.5}
  $$
in the later sum does not contain the factors coming from 
$(p-1)$-st entry.  
We denote the two sum in (3.4.3) by  $\sum_1$  and  $\sum_2$, respectively, 
and denote as
  $$
  =\ {\sum}_1+{\sum}_2\ \ 
  =\ \sum_{w(U)=n}
   \hskip -7pt S_1(U)
  \ \ + \ \ \sum_{r_0=1}^a
   \hskip -20pt
   \sum_{\Sb \hskip 20pt w(U)=n+r_0(p-1) \\ U_{p-1}=0 \endSb}
   \hskip -20pt S_2(U).  
  \tag{3.4.6}
  $$
About  $S_1(U)$  for  $U$  such that  $w(U)=n$, we have
  $$
  \aligned
  w(U[r])+d(U[r])-2
  &=n+(p-1)r+d(U)+r-2 \\
  &=n+pr+d-2 \\
  &=(r+U_{p-1})p+n-pU_{p-1}+d(U)-2. 
  \endaligned
  \tag{3.4.7}
  $$
Note that 
  $$
  \aligned
  n-pU_{p-1}+d(U)-2
 &=(n-(p-1)U_{p-1})+(d(U)-U_{p-1})-2 \\
 &=w(U|_{p-1})+d(U|_{p-1})-2. 
  \endaligned
  \tag{3.4.8}
  $$
About  $S_2(U)$  for  $U$  such that  $w(U)=n+r_0(p-1)$, we have
  $$
  \aligned
  w(U[r-r_0])+&d(U[r-r_0])-2 \\
  =&n+r_0(p-1)+(p-1)(r-r_0)+d(U)+(r-r_0)-2 \\
  =&(r-r_0)p+n+r_0p+d(U)-r_0-2. 
  \endaligned
  \tag{3.4.9}
  $$
According to  $n-pU_{p-1}+d-2$  
(resp. $n+r_0p+d(U)-r_0-2$)  is  $\geqq ap$  or   $< ap$, 
we divide the sum  $\sum_1$  (resp.  $\sum_2$)  into two kinds of sums, 
and we say  $\sum_1=\sum'_1+\sum''_1$  (resp. $\sum_2=\sum'_2+\sum''_2$).  
Here we note that  $n-pU_{p-1}+d(U)-2>0$.  
\fp
(a)  About  $\sum'_1$,  
since  $n-pU_{p-1}+d-2\geqq ap$, 3.2.1 and 3.3.1 yield that  
  $$
  \aligned
  \ord_p(S_1(U))
 &\geqq -\ord_p(\gamma_{\sssize U|_p})+\ord_p((n-pU_{p-1}+d-2)!) \\
 &\geqq \left\lfloor \frac{n-pU_{p-1}+d-2}{2p}\right\rfloor 
  \geqq \left\lfloor \frac{ap}{2p}\right\rfloor 
  =\left\lfloor \frac{a}{2}\right\rfloor. 
  \endaligned
  \tag{3.4.10}
  $$
\fp
(b)  About  $\sum''_1$, we see  $n-pU_{p-1}+d-2< ap$.  
we denote  $N=n-pU_{p-1}+d-2$,  and  $N=pb+e$  ($0\leqq e<p$).  
Lemma 3.3.1 shows that 
  $$
  \aligned
   \ord_p&(N!)
  -\ord_p(\gamma_{\sssize U|_{p-1}}) \\
 &=\ord_p((w(U|_{p-1})+d(U|_{p-1})-2)!)
  -\ord_p(\gamma_{\sssize U|_{p-1}}) \\
 &=\ord_p(\tau_{\sssize U|_{p-1}}) \\
 &\geqq \lfloor N/(2p) \rfloor. 
  \endaligned
  \tag{3.4.11}
  $$
Since  $b< a$, by 3.2.1 we have
  $$
  \aligned
  \ord_p(S_1(U))
  &\geqq 
     a-\Big\lfloor \frac{N}{p} \Big\rfloor
    +\ord_p(N!) 
    -\Big\lfloor\frac{a-\lfloor N/p\rfloor}{p}\Big\rfloor 
    -\ord_p(\gamma_{\sssize U|_{p-1}})\\
  &=\bigg\lfloor \frac{N}{2p} \bigg\rfloor 
  +a-\Big\lfloor \frac{N}{p} \Big\rfloor
    -\Big\lfloor\frac{a-\lfloor N/p\rfloor}{p}\Big\rfloor \\
  &\geqq\Big\lfloor \frac{b}{2} \Big\rfloor+a-b
   -\Big\lfloor \frac{a-b}{p}\Big\rfloor \\
  &> \Big(\frac{b}{2}-1\Big)-b+a-\frac{a-b}{p} \\
  &=-1+\frac{a}{2}+\frac{(a-b)(p-2)}{2p} \\
  &>-1+\frac{a}{2}．  
  \endaligned
  \tag{3.4.12}
  $$
Here the initial side is an integer, 
we have shown that  $\ord_p(S_1(U))\geqq \lfloor a/2 \rfloor$. 
Hence  $\ord_p(\sum_1)\geqq \lfloor a/2 \rfloor$. 
We can prove  $\ord_p(\sum_2)\geqq \lfloor a/2 \rfloor$  
by using 3.2.8 instead of 3.2.1.  
However, we shuld be careful 
for the case  $n+r_0p+d(U)-r_0-2\leqq ap$.  
In this case it should be  $n+r_0p+d(U)-r_0-2\geqq r_0p$  
in order to applying 3.2.8.  
Therefore it must be  $n+r_0p+1-r_0-2\geqq r_0p$  because  $d(U)\geqq 1$.  
Namely, $n-r_0-1\geqq 0$.  
This condition is satisfied for  $r_0=1$, $\cdots$, $a$.  
by our assumption  $n>a$.  
\qed
\vskip 7pt
\fp
{\bf Remark 3.4.13.}\ 
If we replace the condition  $n>a$  in 3.3.1  by  $n\leqq a$,  
we have the similar congruence modulo  $p^{n-1}$  
for the generalized Bernoulli-Hurwitz numbers in \cite{{\bf \^O}}.  
Since those numbers are obtained by specializing 
the universal Bernoulli numbers, the condition  $n>a$  is crucial in 3.3.1.  
\vskip 13pt


\fp
{\bf 3.5. Kummer-Adelberg congruence.} 
\fp
We prove Adelberg's Theorem 3.2 in \cite{{\bf A2}}, 
namely Corollary 3.1.3 of this paper, 
directly from 3.1.1. 
\vskip 3pt
\fp
{\it Proof of {\rm 3.1.3}}. \ 
If  $p=3$, then the statement is vacuous 
by the assumption  $n\not\equiv 0$, $1$ mod $(p-1)$.  
So we may suppose   $p\geqq 5$.  
We prove the desired congruence by induction on $a$.  
The case of  $a=1$  is mentioned in 3.1.2 (1).    
For a given  $a>1$, by taking  $p^{a-1}$  as  $a$  in 3.1.1, we have
  $$
  \sum_{r=0}^{p^{a-1}}
  (-1)^r\binom{p^{a-1}}{r\ \ \ }{A_p}^{p^{a-1}-r}
  \frac{\hat{B}_{n+r(p-1)}}{n+r(p-1)}
  \equiv 0 \mod{p^{\lfloor p^{a-1}/2 \rfloor}}. 
  \tag{3.5.1}
  $$
Because  $a\geqq 2$  and  $p\geqq 5$,  we see  $\lfloor p^{a-1}/2 \rfloor \geqq a$.  
If  $r\neq 0$, $p^{a-1}$, then by (1.1.2) 
  $$
  \aligned
  \ord_p\binom{p^{a-1}}{r\ \ \ }
  &=\frac{S_p(p^{a-1}-r)+S_p(r)-S_p(p^{a-1})}{p-1} \\
  &=\frac{S_p(p^{a-1}-r)+S_p(r)-1}{p-1}. 
  \endaligned
  \tag{3.5.2}
  $$
Let  $\nu=\ord_p(r)$.  By expanding 
  $p^{a-1}-r=d_{a-2}p^{a-2}+d_{a-3}p^{a-3}+\cdots+d_1p+d_0$, 
 ($0\leqq d_j\leqq p-1$)  and  
  $r=h_{a-2}p^{a-2}+h_{a-3}p^{a-3}+\cdots+h_1p+h_0$, 
 ($0\leqq h_j\leqq p-1$) 
$p$-adically, we see obviously that 
  $$
  d_j+h_j=
  \cases
  p-1 & (a-2\geqq j \geqq \nu+1),\\
  p   & (j=\nu), \\
  0   & (\nu-1\geqq j \geqq 0).  
  \endcases
  \tag{3.5.3}
  $$
Hence  $S_p(p^{a-1}-r)+S_p(r)=(p-1)(a-2-\nu)+p$.  
Therefore  
$\ord_p\binom{p^{a-1}}{r\ \ \ }=\ord_p\binom{p^{a-1}}{p^{a-1}-r}=a-1-\nu$.  
Thanks to  $p$  is odd number, we consider the sum 
  $$
  \aligned
   (-1)^r     &   \binom{p^{a-1}}{ r\ \ \ }{A_p}^{p^{a-1}-r}
   \frac{\hat{B}_{n+r(p-1)}}{n+r(p-1)} \\
 &+(-1)^{p^{a-r}} \binom{p^{a-1}}{p^{a-1}-r}{A_p}^{r}
   \frac{\hat{B}_{n+(p^{a-1}-r)(p-1)}}{n+(p^{a-1}-r)(p-1)}
  \endaligned
  \tag{3.5.4}
  $$
for  $1\leqq r \leqq (p^{a-1}-1)/2$.  
Then by the argument above and 2.3.1, we see that the denominator of 
$\frac{\hat{B}_{n+r(p-1)}}{n+r(p-1)}$  for  $0< r < p^{a-1}$  
is not divisible by  $p$.  
Additionally considering the hypothesis of induction, 
we see the sum (3.5.4) is divisible by  $p^a$.  
Thus we have proved 3.1.3.  
\qed



\enddocument
\bye